\documentclass[12pt]{article}
\usepackage{amsthm,amsfonts,amssymb,amscd}

\begin{document}
\begin{center}
\Large \bf Birationally rigid Fano cyclic covers
\end{center}
\vspace{1cm}

\centerline{\large \bf Aleksandr V. Pukhlikov}

\parshape=1
3cm 10cm \noindent {\small \quad \quad \quad
\quad\quad\quad\quad\quad\quad\quad {\bf }\newline We prove
birational superrigidity of Fano cyclic covers over hypersurfaces
in the projective space. } \vspace{1cm}

{\large\bf 1. The main result}
\vspace{0.1cm}

Let ${\mathbb P}={\mathbb P}^{M+1}$ be the complex projective
space, $M\geq 5$, $Q=Q_m\subset{\mathbb P}$ a smooth hypersurface
of degree $m\geq 1$,
$$
\sigma\colon V\to Q
$$
the $K\colon 1$ cyclic cover, branched over a smooth divisor
$W\cap Q$, where $W=W_{Kl}\subset{\mathbb P}$ is a hypersurface of
degree $Kl$. Introducing a new coordinate $u$ of weight $l$, we
can realize $V$ as a complete intersection of type $m\cdot Kl$ in
the weighted projective space
$$
{\mathbb P}(\underbrace{1,\dots,1}\limits_{M+2},l),
$$
namely, $V$ is given by the following system of equations
\begin{equation}\label{1a}
\left\{\begin{array}{l} f(x_0,\dots,x_{M+1})=0\\
u^K=g(x_0,\dots,x_{M+1}),
\end{array}\right.
\end{equation}
where $f(x_*)$ and $g(x_*)$ are homogeneous polynomials of degrees
$m$ and $Kl$, respectively.

If the integers $m$,\,$l$ and $K$ satisfy the relation
$$
m+(K-1)\,l=M+1,
$$
then $V$ is a primitive Fano variety of dimension $M$, that is,
$\mathop{\rm Pic}V={\mathbb Z}K_V$ and $(-K_V)$ is ample.

 The purpose of this note is to sketch a proof of the following

 {\bf Theorem 1.} {\it A general (in the sense of Zariski
 topology) variety $V$ is birationally superrigid.

 In particular, $V$ admits no non-trivial structures of a
 rationally connected fibration, any birational map
 $V\dashrightarrow V^{\sharp}$ onto a Fano variety with ${\mathbb
 Q}$-factorial terminal singularities and $\mathop{\rm
 rk}\mathop{\rm Pic}V^{\sharp}=1$ is an isomorphism and the groups
 of birational and biregular self-maps coincide:}
 $$
 \mathop{\rm Bir}V=\mathop{\rm Aut}V.
 $$
 {\bf Remark 1.} Birational superrigidity means that for any
 movable linear system $\Sigma$ on $V$ its virtual and actual
 thresholds of canonical adjunction coincide:
\begin{equation}\label{1}
c_{\mathop{\rm virt}}(\Sigma)=c(\Sigma),
\end{equation}
see [1-3] for the definitions, also for a (very simple) proof that
(\ref{1}) implies all the other claims of Theorem 1.\vspace{0.3cm}

{\large \bf 2. The regularity conditions} \vspace{0.1cm}

Let us give a precise description of the concept of the variety
$V$ being general. Since $V$ is determined by two polynomials $f$
and $g$, of degrees $m$ and $Kl$, respectively, we can look at $V$
as a point in the parameter space,
$$
(f,g)\in H^0({\mathbb P},{\cal O}_{\mathbb P}(m))\times
H^0({\mathbb P},{\cal O}_{\mathbb P}(Kl))={\cal V}.
$$
Now we will describe the set
$$
U_{\mathop{\rm reg}}\subset {\cal V}
$$
of {\it regular} varieties  $V$, such that $(f,g)\in
U_{\mathop{\rm reg}}$ implies birational superrigidity. The set
$U_{\mathop{\rm reg}}$ is defined by the {\it regularity
conditions} which must be satisfied at every point $o\in V$. These
regularity conditions are similar to those used in [2-5]. We
consider the two cases, when $\sigma(o)\not\in W$ and
$\sigma(o)\in W$, separately.\vspace{0.1cm}

{\bf (R1) The regularity condition outside the ramification
divisor.}

Set $p=\sigma(o)\not\in W$, take a system of affine coordinates
$z_1,\dots,z_{M+1}$ with the origin at $p$ (we can assume that
$z_i=x_i/x_0$), and set $y=u/x^l_0$. Now the standard affine set
$$
{\mathbb A}^{M+2}_{(z_1,\dots,z_{M+1},y)}
$$
is a chart for ${\mathbb P}(1,\dots,1,l)$. With respect to these
coordinates we have two equations of $V$ (abusing our notations,
we denote the corresponding non-homogeneous polynomial by the same
symbol as the homogeneous one):
$$
\begin{array}{l}
f=q_1+\dots+q_m=0,\\
y^K=g=w_0+\dots+w_{Kl},
\end{array}
$$
$q_j$ and $w_J$ are homogeneous polynomials of degree $j$ in
$z_*$.

By assumption, $w_0\neq 0$ and we can assume that $w_0=1$ and
$y(o)=1$. Similar to [2-5], set
$$
g^{1/K}=(1+w_1+\dots+w_{Kl})^{1/K}= 1+\sum^{\infty}_{i=1}\gamma_i
(w_1+\dots+w_{Kl})^{i}=
$$
$$
=1+ \sum^{\infty}_{i=1}\Phi_i(w_1,\dots,w_{Kl}),
$$
where $\gamma_i\in {\mathbb Q}$ are defined by the Taylor
expansion at zero
$$
(1+s)^{1/K}= 1+\sum^{\infty}_{i=1}\gamma_i s^i
$$
and $\Phi_i(w_1(z_*),\dots,w_{Kl}(z_*))$ are homogeneous
polynomials of degree $i\geq 1$ in $z_*$.

Set also
$$
[g^{1/K}]_k=1+ \sum^k_{i=1}\Phi_i(w_1,\dots,w_{Kl}),
$$
and for $k=1,\dots,m$
$$
f_k=q_1+\dots+q_k.
$$

Now we can formulate the regularity condition at the point $o$: if
$m\leq Kl$, then the sequence
$$
q_1,\dots,q_m,\Phi_{l+1}(w_*(z_*)),\dots,\Phi_{Kl-1}(w_*(z_*))
$$
is regular in ${\cal O}_{o,{\mathbb C}^{M+1}}$.

If $m\geq Kl+1$, then we require that the sequence
$$
q_1,\dots,q_{m-1},\Phi_{l+1}(w_*(z_*)),\dots,\Phi_{Kl}(w_*(z_*))
$$
be regular in ${\cal O}_{o,{\mathbb C}^{M+1}}$.\vspace{0.1cm}

{\bf (R2) The regularity condition on the ramification divisor.}

In the notations above, $w_0=0$ in this case. We require that the
sequence
$$
q_1,\dots,q_m,w_1,\dots,w_K
$$
be regular in ${\cal O}_{o,{\mathbb C}^{M+1}}$.

{\bf Proposition 1.} {\it The set $U_{\mathop{\rm reg}}$ of
regular pairs $(f,g)\in{\cal V}$ (that is, the pairs satisfying
(R1) or (R2) at every point $o\in V$) is a non-empty Zariski open
subset in {\cal V}}.

{\bf Proof} is similar to the proof of Proposition 1.1 in [2], see
[5] for more details. \vspace{0.3cm}

{\large \bf 3. The technique of hypertangent divisors}

Here we give a sketch of the proof of Theorem 1.

By the Lefschetz theorem,
$$
\mathbb {\rm Pic}V=A^1V={\mathbb Z}H,
$$
where $H=-K_V$ is the anticanonical class, and
$$
A^2V={\mathbb Z}K^2_V={\mathbb Z}H^2
$$
(the group of codimension 2 cycles modulo numerical equivalence).
By the sufficient condition of birational superrigity (see [2],
Proposition 2.1, also [5,6]), in order to prove Theorem 1 we need
to show that for any irreducible subvariety $Y\subset V$ of
codimension 2 and any point $o\in Y$ the following inequality
\begin{equation}\label{2}
\frac{\mathop{\rm mult}_o}{\mathop{\rm deg}}Y=\frac{\mathop{\rm
mult}_oY}{\mathop{\rm deg}Y}\leq\frac{4}{\mathop{\rm deg}V}
\end{equation}
holds, where $\mathop{\rm deg}Y=(Y\cdot H^{M-2})$ and $\mathop{\rm
deg}V=H^M=mK$.

The inequality (\ref{2}) is proved by means of the technique of
hypertangent divisors [1-6]. Assume at first that the point $o$
lies outside the ramification divisor, $\sigma(o)\not\in W$. For
any $i\leq\mathop{\rm min}(m-1,Kl-1)$ set
$$
\Lambda=|\,\sum^i_{j=0}s_{i-j}f_j+\sum^i_{k=K}s^*_{i-k}(y-[g^{1/K}]_k)\,|
$$
to be the {\it $i$-th hypertangent linear system} at $o$, where
$s_a$, $s^*_a$ are arbitrary homogeneous polynomials in $z_*$ of
degree $a$ and we assume that the value of the coordinate function
$y$ at $o$ is 1 (otherwise replace $y$ by $y\xi$, where
$\xi\in{\mathbb C}^*$ is the appropriate root of 1) and we also
assume that $\sum^i_{k=K}$ is equal to zero if $i\leq K-1$. It is
easy to see that
$$
\Lambda_i\subset |\,iH\,|,\quad\mathop{\rm mult}\nolimits_oD\geq
i+1
$$
for any divisor $D\in \Lambda_i$.

Now we proceed as in [4]: set
$$
{\cal M}=\{1,\dots,m-1\},\quad {\cal L}=\{l,\dots,Kl-2\}.
$$
Here we consider the case $l\geq 3$ and $m\leq Kl$. Set
$$
c_e=\sharp [3,e]\cap{\cal M}+ \sharp [3,e]\cap {\cal L}.
$$
For $e\leq 2$ we get $c_e=0$, for $e\geq \max \{m-1,Kl-2\}$ we get
that $c_e=M-3$. Obviously, $ c_{e+1}\geq c_e.$ Define the {\it
ordering function}
$$
\chi\colon \{1,\dots,M-3\}\to {\mathbb Z}_+
$$
by the formula
\begin{equation}
\label{3} \chi([c_{e-1}+1,c_e]\cap {\mathbb Z}_+)=e.
\end{equation}
If $c_{e-1}=c_e$, then the set $[c_{e-1}+1,c_e]$ is empty and the
formula (\ref{3}) gives no information. Note that
$$
c_{e+1}-c_e\in \{0,1,2\}
$$
so that $\chi$ can take the same value at at most two neighbor
points. Let
$$
{\mathbb D}=\{D_i\in \Lambda_{\chi(i)},\,\,\,i=1,\dots,M-3\}
$$
be a general set of hypertangent divisors.

{\bf Lemma 1.} {\it For every $i=1,\dots,M-3$ (and a sufficiently
general ${\mathbb D}$) the closed algebraic set
$$
R_i({\mathbb D})=\mathop{\bigcap}\limits^i_{j=1} D_i\cap Y
$$
is of codimension $(i+2)$ near the point $o$.}

{\bf Proof:} this follows from the regularity condition (R1), see
[1-6].

Now we construct in the usual way a sequence or irreducible
subvarieties
$$
Y_2=Y,Y_3,\dots,Y_{M-1},
$$
such that
\begin{itemize}

\item $\mathop{\rm codim}Y_i=i$,

\item $Y_{i+1}$ is an irreducible component of the effective
cycle ($Y_i\circ D_{i-1}$) of scheme-theoretic intersection of
$Y_i$ and $D_{i-1}$,

\item the crucial inequality
$$
\frac{\mathop{\rm mult}_o}{\mathop{\rm
deg}}Y_{i+1}\geq\frac{\chi(i-1)+1}{\chi(i-1)}\cdot\frac{\mathop{\rm
mult}_o}{\mathop{\rm deg}}Y_i
$$
\end{itemize}
holds.

Thus
$$
\frac{\mathop{\rm mult}_o}{\mathop{\rm
deg}}Y\cdot\frac{4}{3}\cdot\dots\cdot\frac{m}{m-1}
\cdot\frac{l+1}{l}\cdot\dots\cdot\frac{Kl-1}{Kl-2}=
$$

$$
=\frac{\mathop{\rm mult}_o}{\mathop{\rm
deg}}Y\cdot\frac{m}{3}\cdot\frac{Kl-1}{l}\leq\frac{\mathop{\rm
mult}_o}{\mathop{\rm deg}}Y_{M-1}\leq 1,
$$
whence we get
$$
\frac{\mathop{\rm mult}_o}{\mathop{\rm
deg}}Y\leq\frac{1}{mK}\left(\frac{3Kl}{Kl-1}\right)<\frac{4}{mK},
$$
as required. Thus the estimate (\ref{2}) is proved outside the
ramification divisor.

The remaining cases ($l=2$ or $m\geq Kl$) are treated in a similar
way.\vspace{0.3cm}

{\large \bf 4. The ramification case}\vspace{0.1cm}

Now let us prove the estimate (2) assuming that the point $o$ lies
on the ramification divisor. The hypersurface $W$ is given by the
equation
$$
w_1+\dots+w_{Kl}
$$
(with respect to a fixed system of coordinates $z_*$ with the
origin at $p=\sigma(o)\in{\mathbb P}$). Set
$$
h_k=w_1+\dots+w_k
$$
to be the truncated polynomial, $k\leq K-1$.

{\bf Lemma 2.} {\it The multiplicity of the divisor
($\sigma^*(h_k\,|\,_Q)=0$) at the point $o$ is at least $k+1$}.

{\bf Proof.} This is obvious since
$$
\sigma^*(h_k\,|\,_Q)=(y^K-\sigma^*((w_{k+1}+\dots+w_{Kl})\,|\,_Q))
$$
$y(o)=0$ and $K\geq k+1$. Q.E.D. for the lemma.

Now arguing as in the case outside the ramification divisor, and
applying the regularity condition (R2), we get:
$$
\frac{\mathop{\rm mult}_o}{\mathop{\rm
deg}}Y\cdot\frac{3}{2}\cdot\dots\cdot
\frac{m}{m-1}\cdot\frac{3}{2}\cdot\dots\cdot\frac{K}{K-1}=
$$
$$
=\frac{\mathop{\rm mult}_o}{\mathop{\rm deg}}Y\cdot
\left(\frac{mK}{4}\right)\leq 1,
$$
whence we immediately get the required estimate (\ref{2}).

This completes our proof of the Theorem 1.\vspace{0.3cm}

{\large \bf 5. Further results and concluding
remarks}\vspace{0.1cm}

Using the techniques of [7] and making the regularity conditions
(R1), (R2) stronger in the same way as in [7], one obtains the
following

{\bf Theorem 2.} {\it A sufficiently general (in the sense of
Zariski topology) variety $V\in{\cal V}$ is divisorially
canonical, that is, for any effective divisor $D\in |-nK_V\,|$ the
pair
$$
(V,\frac{1}{n}D)
$$
has canonical singularities.}

In particular, the property of being birationally superrigid is
stable with respect to the operation of taking the direct product
by $V$.

Using the technique of hypertangent divisors (Sec. 3 and 4 above)
in a more delicate way, one can drop the condition of $V$ being a
{\it cyclic} cover of $Q$. Namely, the second equation in the
system (\ref{1a}) can be replaced by the polynomial
$$
u^K+g_1(x_*)u^{K-1}+\dots+g_K(x_*),
$$
where $g_i(x_*)$ are homogeneous of degree $il$.

For $K=2$ this generalization gives nothing new; however, for
$K\geq 3$ it makes the class of varieties much bigger. Denote by
${\cal V}^+$ the corresponding parameter space, ${\cal V}\subset
{\cal V}^+$.

{\bf Theorem 3.} {\it A sufficiently general (in the sense of
Zariski topology) variety $V\in{\cal V}^+$ is divisorially
canonical (in particular, birationally superrigid).}

{\bf Remarks.} (i) The estimate (\ref{2}) can be sharpened at the
expense of making the regularity conditions (R1), (R2) stronger.
This sharpening makes it possible to prove birational rigidity of
Fano fiber spaces $X/{\mathbb P}^1$, the fiber of which is a Fano
multiple hypersurface, see [3,4].

(ii) The standard scheme of proving birational superrigidity is to
establish that for any moving linear system $\Sigma \subset
|-nK_V\,|$ and a {\it general} divisor $D\in \Sigma$ the pair
$$
\left(V,\frac{1}{n}D\right)
$$
has canonical singularities. If this is not the case, then there
exist a birational morphism
$$
\varphi\colon \widetilde V\to V
$$
(one can assume $\varphi$ to be a sequence of blow ups) and an
exceptional divisor $E\subset\widetilde V$ such that the {\it
Noether-Fano inequality}
\begin{equation}\label{4}
\nu_E(\Sigma)> n\cdot a(E)
\end{equation}
is satisfied (such an exceptional divisor is called a {\it maximal
singularity} of the system $\Sigma$). The idea underlying the
sufficient condition of birational superrigidity, used in [1-6]
and in the present paper above, is that if the estimate (\ref{2})
holds, then the Noether-Fano inequality (\ref{4}) is not possible.

It was noted by Cheltsov [8], that if we consider, instead of
rationally connected fiber spaces, the structures of a $K$-{\it
trivial fibration} on $V$, then what we need is just to replace
the Noether-Fano inequality (\ref{4}) by its non-strict version:
\begin{equation}\label{5}
\nu_E(\Sigma)\geq n\cdot a(E).
\end{equation}
But the point is, if the centre
$$
B=\mathop{\rm centre}(E)=\varphi(E)\subset V
$$
of the discrete valuation $E$ on $V$ is a subvariety of
codimension at least 3, then for general divisors $D_1,D_2\in
\Sigma$ the non-strict Noether-Fano inequality (\ref{5}) still
implies that
$$
\mathop{\rm mult}\nolimits_B(D_1\circ D_2)\geq 4n^2
$$
(no modification of the proof is necessary) and the equality takes
place only in the case
\begin{equation}\label{6}
\mathop{\rm codim}B=3,\,\,\,\mathop{\rm mult}\nolimits_B\Sigma=2n,
\end{equation}
so that the estimate (\ref{2}) still gives a contradiction in all
cases but (\ref{6}) and the case $\mathop{\rm codim}B=2$.

Thus, provided that the estimate (\ref{2}) holds for any
irreducible subvariety of codimension two, one has just to study
these two simple cases. For instance, it is easy to see that the
situation (\ref{6}) is impossible for Fano cyclic covers,
considered above. Thus the inequality (\ref{5}) implies that
$B\sim H^2$ and one can show easily that in fact
$$
B=\sigma^{-1}(H_1\cap H_2),
$$
where $H_i$ are hyperplane sections of $Q$, so that all
$K$-trivial structures on $V$ are given by pencils in the
anticanonical linear system $|-K_V\,|$.

This argument is absolutely typical: whenever one has a proof of
birational (super)rigidity, it remains to do a very easy job to
cover the $K$-trivial (or elliptic) case, either. For this reason,
I believe that the $K$-trivial results have very little
independent value. The same applies to the ``Fano structures with
canonical singularities'', see for instance [9] and other papers
in that series.

(iii) As always, one of the crucial parameters of Fano variety is
its anticanonical degree. This can be seen directly from the
inequality (\ref{2}): as the degree gets higher, the required
estimate (2) gets sharper, and, accordingly, harder to prove. On
the other hand, it means that the varieties of small degree are
much easier to investigate. More precisely, primitive Fano
varieties $X$ of degree less or equal than 4 create no problems at
all; in fact, they are all covered (even if they are singular,
provided that the singularities are sufficiently mild) by the test
class construction of V.A.Iskovskikh and Yu.I.Manin in its
original form [10] (plus its higher-dimensional extension [11] or
inversion of adjunction [12]). For this reason, the papers like
[13] are rather of exercise-doing level and hardly represent any
real progress in the field.

(iv) The results of this paper remain true if we take $V$ to be an
iterated Fano cover in the spirit of [5].

(v) One can replace, with no damage to the results, the
hypersurface $Q$ by a complete intersection of type
$d_1\cdot\dots\cdot d_k$, $d_1+\dots+d_k+(K-1)l= M+1$, where $k$
is less than $\frac12\mathop{\rm dim}V$, see [5] and [6].

The detailed proofs of all results discussed above will be
published elsewhere. \vspace{1cm}

\centerline{\bf References} \vspace{0.3cm}

\noindent 1. Pukhlikov A.V., Birational automorphisms of Fano
hypersurfaces, Invent. Math. {\bf 134} (1998), no. 2, 401-426.
\vspace{0.3cm}

\noindent 2. Pukhlikov A.V., Birationally rigid Fano double
hypersurfaces, Sbornik: Mathematics {\bf 191} (2000), No. 6,
101-126. \vspace{0.3cm}

\noindent 3. Pukhlikov A.V., Birationally rigid varieties with a
pencil of Fano double covers. II. Sbornik: Mathematics {\bf 195},
No. 11 (2004), arXiv:math.AG/0403211  \vspace{0.3cm}

\noindent 4. Pukhlikov A.V., Birationally rigid varieties with a
pencil of Fano double covers. I. Sbornik: Mathematics {\bf 195},
No. 7 (2004), arXiv:math.AG/0310270   \vspace{0.3cm}

\noindent 5. Pukhlikov A.V., Birationally rigid iterated Fano
double covers. Izvestiya: Mathematics. {\bf 67} (2003), no. 3,
555-596, arXiv:math.AG/0310268 \vspace{0.3cm}

\noindent 6. Pukhlikov A.V., Birationally rigid Fano complete
intersections, Crelle J. f\" ur die reine und angew. Math. {\bf
541} (2001), 55-79. \vspace{0.3cm}

\noindent 7. Pukhlikov A.V., Birational geometry of Fano direct
products, arXiv:math.AG/0405011.\vspace{0.3cm}

\noindent 8. Cheltsov I.A., On the structures of K-trivial
fibrations on uniruled varieties, Ph.D. Thesis, Steklov Institute,
Moscow and Johns Hopkins University, Baltimore, 1999.
\vspace{0.3cm}

\noindent 9. Cheltsov I.A., Birationally rigid del Pezzo
fibrations, preprint. \vspace{0.3cm}

\noindent 10. Iskovskikh V.A. and Manin Yu.I., Three-dimensional
quartics and counterexamples to the L\" uroth problem, Math. USSR
Sb. {\bf 86} (1971), no. 1, 140-166. \vspace{0.3cm}

\noindent 11. Pukhlikov A.V., Birational automorphisms of a double
space and a double quadric. Math. USSR Izv. {\bf 32} (1989),
233-243. \vspace{0.3cm}

\noindent 12. Koll{\'a}r J., et al., Flips and Abundance for
Algebraic Threefolds, Asterisque 211, 1993. \vspace{0.3cm}

\noindent 13. Cheltsov I.A., Birationally superrigid cyclic triple
spaces, arXiv:math.AG/0410558

\begin{flushleft}
{\it e-mail}: pukh@liv.ac.uk, pukh@mi.ras.ru
\end{flushleft}
\end{document}